\documentclass[12pt]{article} 
\usepackage{amsmath,amssymb,amsthm} 
\usepackage{graphics,epsfig,calc} 
\usepackage[dvipsnames]{color} 
\definecolor{Red}{named}{Red} 
\usepackage{amsmath} 
\usepackage{upgreek} 
\usepackage{latexsym,epsfig,bm,amssymb} 
\usepackage{color} 
\usepackage{amsthm,mathrsfs} 
 
\usepackage{mathptmx} 
 
\newcommand{\beqn}{\begin{eqnarray}} 
\newcommand{\eeqn}{\end{eqnarray}} 
 
 \newcommand{\ba}{\begin{array}} 
\newcommand{\ea}{\end{array}}

\newcommand{\be}{\begin{equation}} 
\newcommand{\ee}{\end{equation}} 

\newcommand{\cE}{{\cal E}} 
\newcommand{\cX}{{\cal X}} 
\newcommand{\si}{\sigma} 
\newcommand{\de}{\delta} 
\newcommand{\la}{\label} 
\newcommand{\cV}{{\cal V}} 
\newcommand{\Lam}{\Lambda} 
\newcommand{\cR}{{\cal R}} 
 \newcommand{\cK}{{\cal K}} 
\newcommand{\re}{\ref} 
\newcommand{\ve}{\varepsilon} 
\newcommand{\cW}{{\cal W}} 
\newcommand\R{{\mathbb R}} 
\newcommand{\ci}{\cite} 
\newcommand{\fr}{\frac} 
\newcommand{\De}{\Delta} 
\newcommand{\ds}{\displaystyle} 
\newcommand{\om}{\omega} 
\newcommand{\cS}{{\cal S}}

\newcommand{\bo}{{\hfill\loota}}

\newcommand{\loota}{\hbox{\enspace{\vrule height 7pt depth 0pt width 7pt}}}

\newcommand{\5}{{\hspace{0.5mm}}}

\newcommand{\dimension}{\operatorname{dim}} 
 
\newcommand{\sgn}{\mathop{\rm sgn}\nolimits}

\newcommand{\rIm}{{\rm Im\5}} 
 
\renewcommand{\Pr}{\hspace{-6mm}{\bf Proof.~}} 
 
\newcommand{\Ker}{{\rm Ker\5}} 
\newcommand{\Ran}{{\rm Ran\5}}

\renewcommand{\theequation}{\thesection.\arabic{equation}} 
\newtheorem{theorem}{Theorem}[section] 
\renewcommand{\thetheorem}{\arabic{section}.\arabic{theorem}} 
\newtheorem{definition}[theorem]{Definition} 
 
\newtheorem{lemma}[theorem]{Lemma} 
\newtheorem{example}[theorem]{Example} 
\newtheorem{remark}[theorem]{Remark} 
\newtheorem{remarks}[theorem]{Remarks} 
\newtheorem{cor}[theorem]{Corollary} 
\newtheorem{pro}[theorem]{Proposition}

\newcommand{\bd}{\begin{definition}} 
 \newcommand{\ed}{\end{definition}} 
\newcommand{\bt}{\begin{theorem}} 
 \newcommand{\et}{\end{theorem}} 
\newcommand{\bqt}{\begin{qtheorem}} 
 \newcommand{\eqt}{\end{qtheorem}}

\newcommand{\bp}{\begin{pro}} 
 \newcommand{\ep}{\end{pro}} 
 
\newcommand{\bl}{\begin{lemma}} 
 \newcommand{\el}{\end{lemma}} 
\newcommand{\bc}{\begin{cor}} 
 \newcommand{\ec}{\end{cor}} 
 
\newcommand{\bex}{\begin{example}} 
 \newcommand{\eex}{\end{example}} 
\newcommand{\bexs}{\begin{examples}} 
 \newcommand{\eexs}{\end{examples}}

\newcommand{\bexe}{\begin{exercice}} 
 \newcommand{\eexe}{\end{exercice}}

\newcommand{\br}{\begin{remark} } 
 \newcommand{\er}{\end{remark}} 
\newcommand{\brs}{\begin{remarks}} 
 \newcommand{\ers}{\end{remarks}} 

\begin{document} 
\begin{titlepage}

\begin{center} 
{\Large\bf 
On the eigenfunction expansion for 
\bigskip\\ 
Hamilton operators} 
\end{center} 
\vspace{1cm} 
 \begin{center} 
{\large A. Komech}  
\\ 
{\it Faculty of Mathematics of Vienna University\\ 
and Institute for Information Transmission Problems RAS } \\ 
e-mail:~alexander.komech@univie.ac.at 
\bigskip\\ 
{\large E. Kopylova} 
\\ 
{\it Faculty of Mathematics of Vienna University\\ 
and Institute for Information Transmission Problems RAS} \\ 
 e-mail:~elena.kopylova@univie.ac.at 
\end{center} 
\vspace{1cm}

 \begin{abstract}

A spectral representation for solutions to linear Hamilton equations 
with nonnegative energy in Hilbert spaces is obtained. 
This paper continues our previous work on Hamilton equations 
with positive definite  energy. 
Our approach is a~special version of M.~Krein's 
spectral theory of $J$-selfadjoint operators in 
Hilbert spaces with indefinite metric. 
 
As a principal application of these results, we justify  
the eigenfunction expansion for 
linearized nonlinear relativistic Ginzburg--\allowbreak Landau equations. 
 
\bigskip\bigskip\bigskip 
 
{\it Key words and phrases}: Hamilton equation; selfadjoint operator; 
$J$-selfadjoint operator; Krein space; spectral resolution; 
spectral representation; secular solutions; eigenvector; Jordan block; 
 Ginzburg--\allowbreak Landau equation; kink; asymptotic stability; 
generalized eigenfunction; eigenfunction expansion; Fermi Golden Rule. 
\end{abstract}

{\it 2010 Mathematical Subject Classification}: 35P, 37K

\end{titlepage}

\section{Introduction} 
 
We consider complex linear Hamilton operators in a complex Hilbert space $\mathcal X$, 
\begin{equation}\label{Hs} 
A=JB,\qquad \mbox{where }\qquad B^*=B, \qquad J^*=-J, 
\qquad J^2=-1. 
\end{equation} 
In particular, the 
operator $J:\cX\to\cX$ is  bounded. 
The selfadjoint operator 
 $B$ is defined on a dense domain $D(B)\subset\mathcal X$. 
Our aim is to prove the well-posedness of the Cauchy problem for the equation 
\begin{equation}\label{NHl} 
\dot X(t)=AX(t), 
\end{equation} 
and obtain a spectral representation for solutions 
and the corresponding spectral resolution for $A$. 
For example, for $J=i$ the solutions are given by 
$X(t)=e^{iBt}X(0)$. 
A~more general  `commutative case', when $JB=BJ$, 
reduces to $J=i$, since $JB=iB_1$, where $B_1=-iJB$ is the selfadjoint operator. 
However, $JB\ne BJ$ for linearizations of 
$U(1)$-invariant 
nonlinear Schr\"odinger equations as shown in Appendix of~\cite{KK13}. 
 
We develop the theory in the case of nonnegative `energy operators' 
$B$ with spectral gap and finite `degeneracy of the vacuum': 
 
\smallskip 
\textbf{Condition~I}\hfill $\sigma (B)\subset \{0\}\cup [\delta,\infty)$, \ \ $\delta>0$\hfill (1.3) 
 
\smallskip 
 
\textbf{Condition~II}\hfill $\dimension\Ker B<\infty$.\hfill (1.4) 
\smallskip 
\setcounter{equation}{4}~\\ 
These conditions hold, in particular, for all equations 
considered in \cite{IKS12}--\cite{KKm11}. 
The motivation for the theory was discussed in \cite{KK13}, 
in which the simplest 
case $\si (B)\subset (\de,\infty)$ (i.e., $\dim~ \Ker B=0$) 
 was studied.

We reduce the problem to a selfadjoint generator 
developing  a~special version of M.~Krein's spectral theory 
of $J$-selfadjoint operators in Hilbert spaces with indefinite metric \cite{AI1989},~\cite{KL1963}. 
We apply this version for justification of the 
eigenfunction expansions for the linearization of 
relativistic nonlinear Ginz\-burg--\allowbreak Landau 
equation \cite{KK11}. The generator of the linearization reads 
\begin{equation}\label{AB0i} 
A=\left( 
\begin{array}{cc} 
0   & 1 \\ 
-S & 0 
\end{array} 
\right), 
\end{equation} 
where $S:=-\frac{d^2}{dx^2}+m^2+ V_0(x)$. 
Our results are concerned with the following: 
 
$\bullet$ 
The existence and uniqueness and formula for generalized solutions 
to \eqref{NHl} under conditions (1.3), (1.4) 
for all initial states~$X$ with finite 
energy $\langle B X,X\rangle$. Here, $\langle \cdot,\cdot \rangle$ stands for 
the scalar product in~$\mathcal X$. 
 
$\bullet$ The eigenfunction expansion 
\begin{equation}\label{uns30i} 
\left( 
\begin{array}{c} 
\psi(t)\\ 
\dot \psi(t) 
\end{array} 
\right) 
=t\Phi_0+\Psi_0+ 
\sum e^{-i\omega_k t} C_k a_k+ 
\int_{|\omega|\ge m} e^{-i\omega t} C(\omega)a_\omega~d\omega 
\end{equation} 
for solutions to \eqref{NHl} with generator (\ref{AB0i}). 
Here,  $\Phi_0\in\Ker A$, and 
$a_k$ are the eigenvectors of $A$, 
$\Psi_0$ is the associated eigenvector to 
$\Phi_0$, while 
$a_\omega$ are 
generalized eigenfunctions of~$A$. 
 
Such eigenfunction expansions were used in \cite{BP95,BS03,KK11} 
for the calculation of `Fermi Golden Rule' (FGR) in the context of the 
nonlinear 
Schr\"odinger and Klein--\allowbreak Gordon equations. 
This is a nondegeneracy condition, which was introduced 
in~\cite{Sigal93} in the framework of nonlinear wave and Schr\"odinger 
equations. 
This condition means a~strong coupling of discrete 
and continuous spectral components of solutions providing the radiation of energy to infinity and 
which results in the asymptotic stability of solitary waves. 
The calculation of FGR, as given in \cite{BP95,BS03,KK13}, relies 
on eigenfunction expansions of type 
(\ref{uns30i}). Our main Theorem \ref{tmain} justifies the eigenfunction expansion 
\cite[(5.14)]{KK11}, for which no detailed proof was given before. 
This justification was one of our main motivation for writing the present paper. 
 
The eigenfunction expansion (\ref{uns30i}) 
extends our previous result  \cite{KK13}, where the expansion was established 
only for odd solutions. In this framework 
we have $\Ker B=0$ 
and $\Phi_0=\Psi_0=0$. This framework was sufficient for the proof of 
asymptotic stability of 
standing solitons for the nonlinear relativistic Ginzburg--\allowbreak Landau equations 
under odd perturbations \cite{KK11}. 
However, to establish the asymptotic stability under 
arbitrary perturbations we need the expansion  (\ref{uns30i}) 
for solutions without antisymmetry. 
 
\medskip 
 
Let us comment on our approach. First, we 
reduce the abstract problem 
\eqref{NHl} under conditions (1.3), (1.4) 
to a selfadjoint generator 
justifying the classical M.~Krein transformation~\cite{GK}. 
This reduction is a~special version of spectral theory 
of $J$-selfadjoint operators in Hilbert spaces with indefinite metric 
 \cite{AI1989,KL1963}, extending our approach \cite{KK13} to the case 
$\Ker B\ne 0$. 
This extension required new robust ideas 
i) to  analyze the structure 
of spectrum of the reduced selfadjoint operator, and 
ii) to find the canonical form of the Hamilton operator. 
We provide a broad range of examples 
 satisfying all the imposed conditions (1.3),(1.4), (2.12), and~(3.1).

Second, we apply this abstract spectral theory 
to operator (\ref{AB0i}) and construct the 
eigenfunction expansion 
for the reduced selfadjoint operator 
following the method of Section 5 from~\cite{KK13}. 
At last, we deduce (\ref{uns30i}) by extending our approach from~\cite{KK13}, which relies on the methods of~PDO. 
 
One of our  novelties is a vector-valued 
treatment of the convergence of the integral over the continuous spectrum in (\ref{uns30i}). 
Namely, we show 
that the integral is the limit of the 
corresponding integrals over $m\le|\omega|\le M$   as $M\to\infty$ 
in the Sobolev space $H^1({\mathbb R})$. 
In its own turn, the integral over $m\le|\omega|\le M$ is absolutely 
converging in the weighted $L^2$-space with the weight $(1+|x|)^{-s}$, 
where $s>1$. 
 
Finally, calculation of the symplectic normalization of the 
generalized eigenfunctions requires extra arguments pertaining to the nondegenerate case~\cite{KK13}. 
 
\medskip

We now give some comments on the related works. 
Some spectral properties of the Hamilton non-selfadjoint 
operators were studied by 
V.~Buslaev and G.~Perelman \cite{BP93,BP95,BS03}, 
M.\,B.~Erdogan and W.~Schlag \cite{ES06,S07}, S.~Cuccagna, D.~Pelinovsky and 
V.~Vougalter \cite{CPV05}. 
It is worth noting that the eigenfunction expansions of $J$-selfadjoint operators were 
not justified previously. 
 
Spectral resolution 
of bounded $J$-selfadjoint nonnegative 
operators in Krein spaces was constructed 
 by M.~Krein, H.~Langer and Yu.~Shmul'yan \cite{KL1963,KS1966}, 
and extended to unbounded {\it definitizable} operators by M.~Krein, 
P.~Jonas, H.~Langer and others 
\cite{ IKL1982,  Jonas1981, Jonas1988,   L1981, LN1983}. 
The corresponding unitary operators were examined by P.~Jonas \cite{Jonas1986}. 
However, the spectral resolution  alone is insufficient 
for justification of eigenfunction expansion. 
Our version of the theory under conditions (1.3), (1.4) 
allows us to justify the eigenfunction expansion (\ref{AB0i}). 
 
The spectral theory of 
definitizable operators was applied to the Klein--\allowbreak Gordon equations 
with non-positive energy 
by P.~Jonas, H.~Langer, B.~Najman and C.~Tretter 
\cite{Jonas1993, Jonas2000, LNT2006,LNT2008,LT2006}, 
where the existence and uniqueness 
of classical solutions were proved, 
and the existence of unstable eigenvalues (imaginary frequencies) 
was studied. The instability is related to the known {\it Klein paradox} 
in quantum mechanics \cite{Sakurai}. 
 
The scattering theory 
for the  Klein--\allowbreak Gordon equations 
with non-positive energy 
was developed by C.~G\'erard and T.~Kako 
using the theory of  definitizable operators in Krein spaces \cite{Gerard2012, Kako1976}. 
 
\medskip 
 
The plan of our paper is as follows. 
In Section~2, we justify the M.~Krein transformation 
under conditions (1.3), (1.4), and find the structure 
of spectrum of the corresponding selfadjoint generator. 
In Section~3,  we construct the spectral representation 
for solutions to \eqref{NHl} and 
deduce the canonical form 
of the Hamilton generator. 
In Section~4, we check all conditions 
 (1.3), (1.4), (2.12), and (3.1) 
for operator~(\ref{AB0i}). 
In Sections 5 and~6, we justify 
the eigenfunction expansion 
(\ref{uns30i}) by applying the methods of Sections 3--4. 
In Section~7, we calculate symplectic normalization 
of the generalized eigenfunctions. 
Finally, in the Appendix  we construct examples of Hamilton equations 
satisfying all the imposed conditions. 
 
\medskip 
 
{\bf Acknowledgments.} The authors take pleasure in thanking A.~Kostenko 
and G.~Teschl for useful discussions on spectral theory of  $J$-selfadjoint operators.

A.K. was supported partly 
by Alexander von Humboldt Research Award, 
Austrian Science Fund (FWF): P22198-N13, 
and the Russian Foundation for Basic Research.

E.K. was supported partly by Austrian Science Fund (FWF): M1329-N13, 
and the Russian Foundation for Basic Research.

\setcounter{equation}{0} 
\section{Reduction to  symmetric generator} 
 
In this section, we shall reduce \eqref{NHl} to an equation with selfadjoint generator. 
 
\subsection{Generalized solutions} 
 
Throughout the paper, $D(B)$ is 
a dense domain of the selfadjoint operator $B$. 
We set $\Lambda:=B^{1/2}\ge  0$ and denote by 
$\mathcal V\subset \mathcal X$ 
the Hilbert space 
 which is the domain of $\Lambda$ 
endowed with the norm 
\begin{equation}\label{Vn22} 
\Vert X\Vert_{\mathcal V} 
:=\Vert \Lambda X\Vert_\mathcal X+\Vert X\Vert_\mathcal X. 
\end{equation} 
We have the continuous injections  of Hilbert spaces 
$\mathcal V\subset \mathcal X$, and 
the operator 
\begin{equation}\label{LL} 
\Lambda: \mathcal V\to \mathcal X 
\end{equation} 
is continuous. By definition (\ref{Vn22}), 
\begin{equation}\label{KV} 
\mathcal K\subset \mathcal V. 
\end{equation} 
For example, $\mathcal V$ becomes the Sobolev space 
$H^1({\mathbb R}^n)$ if $\mathcal X=L^2({\mathbb R}^n)$ and 
$A= -i\Delta$. 
 
Since $\Lambda$ and $B$ are selfadjoint operators, 
we have 
\begin{equation}\label{KR} 
\mathcal X=\mathcal K\oplus \mathcal R,~~
\mathcal K:=\Ker \Lambda=\Ker B,~~  \mathcal R:=\Ran \Lambda=\Ran B 
= 
\mathcal K^\bot. 
\end{equation} 
Further, we assume henceforth that  $\mathcal R$ is endowed with the norm of $\mathcal X$. 
Then 
 $\Lambda_+:=\Lambda |_\mathcal R: \mathcal R\cap\mathcal V \to \mathcal R$ 
is an invertible operator by (1.3); i.e., 
\begin{equation}\label{LR} 
\Lambda_+^{-1}: \mathcal R \to \mathcal V 
\end{equation} 
is the bounded operator. 
We will consider solutions 
\begin{equation}\label{YD} 
X(t)\in C({\mathbb R},\mathcal V) 
\end{equation} 
to equation \eqref{NHl}. 
The equation will be understood in the sense of {\it mild solutions} 
\cite{Cazenave} 
\begin{equation}\label{ini} 
X(t)-X(0)=A\int_0^t X(s)ds,\qquad t\in{\mathbb R}, 
\end{equation} 
where the Riemann integral converges in $\mathcal V$ by  (\ref{YD}). 
 
\subsection{Krein substitution} 
Let us reduce equation 
\eqref{ini} by the well-known substitution 
\begin{equation}\label{ZY} 
Z(t):=\Lambda X(t)\in C({\mathbb R},\mathcal R) 
\end{equation} 
used by M. Krein in 
the theory of parametric resonance: 
see formula  (1.40) of 
\cite[Chapter VI]{GK}. 
Applying  $\Lambda$ to both sides of 
equation \eqref{ini}, we obtain 
\begin{equation}\label{L2} 
Z(t)-Z(0)=\Lambda J\Lambda\int_0^t Z(s)ds, \qquad t\in{\mathbb R}. 
\end{equation} 
Formally, 
\eqref{L2} reads 
\begin{equation}\label{CPF3} 
  i\dot Z(t)=H 
Z(t), \qquad t\in{\mathbb R}, 
\end{equation} 
where $H$ stands for the `Schr\"odinger operator' 
\begin{equation}\label{CH} 
H=\Lambda iJ\Lambda, 
\end{equation} 
which is 'formally symmetric'. 
 
\subsection{Equivalence of reduction} 
 
In order to prove the equivalence of equations  \eqref{NHl} and \eqref{CPF3} 
we introduce the following new condition. 
 
\smallskip 
\textbf{Condition~III}\hfill $J\mathcal K\subset \mathcal V$.\hfill (2.12) 
\smallskip
\setcounter{equation}{12}~\\
We denote by $\Pi_\mathcal K:\mathcal X\to\mathcal K$ 
 the orthogonal projection, and set 
$$ 
P:=\Pi_\mathcal K J\Lambda. 
$$ 
\bl\label{lP} 
Let conditions  {\rm (1.4)} and  {\rm (2.12)} hold. Then the operator 
$P: \mathcal X\to\mathcal V$ is continuous. 
\el 
 \Pr 
It suffices to note that 
\begin{equation}\label{cont} 
\Pi_\mathcal K J\Lambda=\sum_1^N|Y_k\rangle\langle Y_k|J\Lambda 
=-\sum_1^N|Y_k\rangle\langle \Lambda J Y_k|, 
\end{equation} 
where $Y_k\in\mathcal K\subset \mathcal V$, 
$N=\dimension\mathcal K$, 
 and $\Lambda J Y_k\in\mathcal X$ by  {\rm (2.12)}. 
\bo 
\medskip 
 
Equation \eqref{ini} 
with  $X(t)\in C({\mathbb R},\mathcal V)$ 
can be written as 
\begin{equation}\label{iniZ} 
X(t)-X(0)=J\Lambda \int_0^t Z(s)ds, \qquad t\in{\mathbb R}. 
\end{equation} 
By Lemma \ref{lP} 
this equation  implies the system 
\begin{equation}\label{inisys} 
\left.
\ba{rcl}
X_\mathcal R(t)-X_\mathcal R(0) 
&=&(1-\Pi_\mathcal K) J\Lambda 
\ds\int_0^t Z(s)ds
\\
\\
\qquad X_\mathcal K(t)-X_\mathcal K(0) 
&=&\Pi_\mathcal K J\Lambda \ds\int_0^t Z(s)ds, 
\ea\right|
\end{equation} 
where $X_\mathcal K(t)=\Pi_\mathcal K X(t)$ and   $X_\mathcal R(t)=(1-\Pi_\mathcal K) X(t)$.

\bl\label{leq1} 
i) Let $X(t)\in C({\mathbb R},\mathcal V)$ be a solution to \eqref{NHl} in the sense \eqref{ini}. 
Then $Z(t)=\Lambda X(t)\in C({\mathbb R},\mathcal R)$ 
is the solution to \eqref{CPF3} in the sense~\eqref{L2}. 
\medskip\\ 
ii) Let  $Z(t)\in C({\mathbb R},\mathcal R)$  be a fixed solution to \eqref{CPF3} in the sense  \eqref{L2}. 
Then  there exists a unique solution 
$X(t)\in C({\mathbb R},\mathcal V)$ to \eqref{NHl} in the sense \eqref{ini} 
satisfying  \eqref{ZY}.

\el 
\Pr 
It suffices to prove ii). The uniqueness holds, because 
\begin{equation}\label{inisys3} 
X_\mathcal R(t) 
=\Lambda_+^{-1} Z(t), \qquad 
X_\mathcal K(t)-X_\mathcal K(0) 
=\Pi_\mathcal K J\Lambda \int_0^t Z(s)ds, 
\end{equation} 
where the first equation follows from 
\eqref{ZY}, and the second one, from the second equation of~(\ref{inisys}). 
 
To prove the existence we define $X_\mathcal R(t)$ and 
$X_\mathcal K(t)$ by (\ref{inisys3}). 
Then \eqref{ZY} holds, and 
$X(t)=X_\mathcal R(t)+X_\mathcal K(t)\in C({\mathbb R},\mathcal V)$. Hence, 
the first equation (\ref{inisys3}) together with 
\eqref{L2} and \eqref{ZY} 
imply that 
\beqn\label{inisys4} 
X_\mathcal R(t)-X_\mathcal R(0) 
&=&\Lambda_+^{-1} [Z(t)-Z(0)]=\Lambda_+^{-1} \Lambda J\Lambda\int_0^t Z(s)ds
\nonumber\\
\nonumber\\
&=& 
(1-\Pi_\mathcal K) J\Lambda^2 
\int_0^t X(s)ds. 
\eeqn 
Finally, the second equation (\ref{inisys3}) can be written as 
\begin{equation}\label{inisys5} 
X_\mathcal K(t)-X_\mathcal K(0) 
=\Pi_\mathcal K J\Lambda^2 \int_0^t X(s)ds~ 
\end{equation} 
by \eqref{ZY}. 
Summing up, we obtain \eqref{ini}. 
\bo

\subsection{Symmetry and spectrum}

The  domain of $H$ is equal to 
\be\la{DH} 
D(H)=\{Z\in\cV: J\Lam Z\in\cV\}=\Lam^{-1}_\cR(J\cV\cap\cR)+\cK~. 
\ee 
Obviously, the operator $H$ is symmetric on  $D(H)$, and hence, 
 $H$ is a closable operator in $\cX$. 
However, 
we still do not know whether its domain is dense in $\cX$. 
This is why we need our last condition 
 
\smallskip 
\textbf{Condition~IV}\hfill $H^*=H$.\hfill (2.21) 
\smallskip 
\setcounter{equation}{21}~\\
A broad range of examples is provided by 
Lemma \re{lcond23}. A concrete example is given by (\re{LH0}).

\bt\la{tspec} 
Let conditions {\rm (1.3),  (1.4),   (2.12)},  and {\rm (2.21)}  hold. 
Then 
\be\la{spec} 
\si(H)\subset 
(-\infty,-\ve]\cup 0 \cup [\ve,\infty) 
\ee 
with some  $\ve>0$. 
\et 
\Pr 
 Operator $\Lam+\Pi_\cK:\cV\to\cX$ is invertible by condition (1.3), 
since $(\Lam+\Pi_\cK)|_\cK=\Pi_\cK$ and $(\Lam+\Pi_\cK)|_\cR=\Lam_+$. 
Hence, the operator 
\begin{equation}\label{H1} 
H_+:=(\Lambda+\Pi_\mathcal K)iJ(\Lambda+\Pi_\mathcal K) 
\end{equation} 
is also invertible; i.e., its inverse 
\begin{equation}\label{H2} 
H_+^{-1}:=(\Lambda+\Pi_\mathcal K)^{-1}iJ(\Lambda+\Pi_\mathcal K)^{-1} 
\end{equation} 
is a~bounded operator on~$\mathcal X$. 
On the other hand, this operator is symmetric on $\mathcal X$, and hence 
it is selfadjoint. Moreover, $H_+$ is injective operator 
on $\mathcal X$. 
Hence, Theorem 13.11 (b) of \cite{Rudin} implies that 
$H_+$ is a~selfadjoint operator with a 
dense domain $D(H_+)$. Further, 
\beqn
\label{H3} 
H_+&=&H+\Pi_\mathcal K iJ(\Lambda+\Pi_\mathcal K)+(\Lambda+\Pi_\mathcal K)iJ\Pi_\mathcal K 
+\Pi_\mathcal K iJ\Pi_\mathcal K
\nonumber\\
\nonumber\\
&=&H+T. 
\eeqn
Here, $\Pi_\mathcal K iJ(\Lambda+\Pi_\mathcal K)$ and 
$\Pi_\mathcal K iJ\Pi_\mathcal K$ 
are finite-range operators $\mathcal V\to\mathcal V$. 
On the other hand, {\rm (2.12)} implies that 
$(\Lambda+\nobreak \Pi_\mathcal K)iJ\Pi_\mathcal K$ is also a~finite-range operator from $\mathcal V$ to~$\mathcal V$. 
Hence, $T:\mathcal V\to\mathcal V$ is the finite-range operator 
which is symmetric in $\mathcal X$. As the result, 
(\ref{H3}) implies that $H$ is defined 
and symmetric 
on  $D(H_+)$. 
 
Further, 
the resolvent 
$(H_+-\lambda)^{-1}:\mathcal X\to \mathcal X$ is bounded and analytic 
 in a small complex neighborhood  $\mathcal O$ 
of $\lambda=0$, and 
\begin{equation}\label{HT} 
H-\lambda=H_+-\lambda-T=[1-T(H_+-\lambda)^{-1}](H_+-\lambda), 
~~~~~\lambda\in\mathcal O. 
\end{equation} 
Here, the operator $H-\lambda$ is invertible for $\rIm\lambda\ne 0$ by {\rm (2.21)}, 
while $H_+-\lambda$ is invertible in a small complex 
neighborhood  $\mathcal O$ 
of $\lambda=0$. 
Hence, $\Ker [1-T(H_+-\lambda)^{-1}]=0$ for $\lambda\in\mathcal O$ 
with $\rIm\lambda\ne 0$. 
Therefore, 
$1-T(H_+-\lambda)^{-1}$ 
is invertible for these~$\lambda$ 
by Fredholm's theorem,  inasmuch as $T$ is a~finite-rank operator. 
Hence, it is also invertible 
in~$\mathcal O$ outside a discrete set. 
Now (\ref{HT}) implies (\ref{spec}).\bo 
 
\br\label{sad} 
Let conditions {\rm (1.3), (1.4)}, and {\rm (2.12)} hold. 
Then 
\medskip\\ 
i) The domain of $H$ is dense in $\mathcal X$, as is shown in the proof of Theorem~{\rm \ref{tspec}}. 
\medskip\\ 
ii) $H$ admits selfadjoint 
extensions, because 
\begin{equation}\label{Npm} 
N_+=N_-, \qquad N_\pm:=\dimension [\Ran (H\mp i)]^\bot. 
\end{equation} 
Indeed, $\Ran(H_+-\lambda)=\mathcal X$ for  $\lambda$ from 
a small complex neighborhood $\mathcal O$ of $\lambda=0$. 
On the other hand, 
the dimension of $({\Ran [1-T(H_+-\lambda)^{-1}]})^\bot$ is constant 
in  $\mathcal O$ 
outside a discrete set, because $T$ is a~finite-rank operator. 
Therefore,  (\ref{HT}) implies that 
$\dimension [\Ran (H-\lambda)]^\bot$ 
is also constant in  $\mathcal O$ 
outside a discrete set, verifying~(\ref{Npm}). \bo 
\medskip

\er

\setcounter{equation}{0} 
\section{Dynamical group and canonical form} 
We construct spectral representation for solutions to   \eqref{NHl} 
and deduce the canonical form of the Hamilton generator. 
 
\subsection{Spectral representation of solutions} 
 
We will construct solutions to \eqref{CPF3}, and afterwords, 
reconstruct the corresponding solutions to \eqref{NHl}. The Spectral Theorem implies the following lemma.

\bl\label{ceqv} 
Let conditions {\rm (1.3), (1.4), {\rm (2.12)}, {\rm (3.1)}} hold. Then, 
for any $Z(0)\in\mathcal R$, equation {\rm \eqref{CPF3}} 
admits a unique solution $Z(t)\in C({\mathbb R},\mathcal R)$ in the sense {\rm \eqref{L2}}. 
The solution is given by 
\begin{equation}\label{uns222} 
Z(t)=e^{-iHt}Z(0)\in C({\mathbb R},\mathcal R). 
\end{equation} 
\el 
 
Now we can reconstruct solutions to \eqref{NHl} using formulas 
(\ref{inisys3}): 
\begin{equation}\label{unsX} 
X(t)=\Lambda_+^{-1}e^{-iHt}\Lambda X(0)  +X_\mathcal K(0) 
+ 
P \int_0^t e^{-iHs} \Lambda X(0)ds, 
\end{equation} 
where the operator 
$P:\mathcal R\to\mathcal V$ is  bounded  by Lemma \ref{lP}. 
To 
evaluate the integral in (\ref{unsX}), 
we denote by $\Pi_0$  and $\Pi_R$, respectively, the spectral projections of~$\mathcal X$ 
onto $\Ker H\cap \mathcal R$ and $R:=\Ran H\subset\mathcal R$. 
Obviously, 
\begin{equation}\label{PI} 
e^{-iHs}=\Pi_0+e^{-iH_Rs}\Pi_R,~~ \int_0^te^{-iHs}ds=t\Pi_0+ 
i(e^{-iH_Rt}-1)H_R^{-1}\Pi_R, 
\end{equation} 
where $H_R:=H|_{R\cap D(H)}$. 
Now (\ref{unsX}) reads 
\beqn\label{unsXr} 
X(t)&=&\Lambda_+^{-1}e^{-iHt}\Lambda X(0)  +X_\mathcal K(0) 
\nonumber\\
\nonumber\\
&&+
t P \Pi_0\Lambda X(0)+ 
i P (e^{-iH_Rt}-1)H_R^{-1}\Pi_R \Lambda X(0). 
\eeqn 
Lemmas  \ref{ceqv} and \ref{leq1} imply the following proposition.

\bp\label{pX} Let conditions {\rm (1.3), (1.4), {\rm (2.12)}, {\rm (3.1)}} hold. 
Then, for any  $X(0)\in\mathcal V$, 
\medskip\\ 
i) Equation \eqref{NHl} 
admits  a unique solution $X(t)\in C({\mathbb R},\mathcal V)$. 
\medskip\\ 
ii) The solution admits 
the spectral representation \eqref{unsXr}. 
\ep

\subsection{Spectral resolution} 
 
Representation \eqref{unsXr} can be written as 
\beqn\label{unsXrw} 
e^{At}&=& 
\Lambda_+^{-1}\int_{\mathbb R} 
e^{-i\omega t}dE(\omega)\Lambda   +\Pi_\mathcal K 
\nonumber\\
\nonumber\\
&&+ 
t P \Pi_0\Lambda + 
i P \int_{|\omega|\ge \varepsilon}\frac{e^{-i\omega t}-1}{\omega} dE(\omega) \Lambda , 
\eeqn 
where $dE(\omega)$ denotes the spectral family of $H$, and $\varepsilon>0$ is the number 
from (\ref{spec}). 
Formally, 
\begin{equation}\label{unsXrw2} 
e^{At}=\int_{\mathbb R} e^{-i\omega t}dF(\omega), 
\end{equation} 
where 
\beqn\label{unsF} 
dF(\omega) 
&=& 
\Big[\Lambda_+^{-1} 
+ 
\frac{iP}{\omega}\Big] 
\chi_\varepsilon(\omega)dE(\omega)\Lambda 
\nonumber\\ 
\nonumber\\ 
&&+ 
\Big[\Pi_\mathcal K- 
i P \int_{|\omega|\ge \varepsilon}\frac 
{dE(\omega)} 
{\omega}  \Lambda 
\Big]\delta(\omega)d\omega 
\nonumber\\ 
\nonumber\\ 
&&  +\Lambda_+^{-1}\Pi_0\Lambda\delta(\omega)d\omega 
-i 
P \Pi_0\Lambda \delta'(\omega)d\omega, 
\eeqn 
and $\chi_\varepsilon$ is the indicator of the set $|\omega|\ge \varepsilon$. 
Setting $t=0$ in 
both sides of 
(\re{unsXrw2}) and in their derivatives, 
we formally obtain 
\begin{equation}\label{unsXrw3} 
1=\int_{\mathbb R} 
dF(\omega),\qquad A=-i\int_{\mathbb R} \omega\5\5 dF(\omega). 
\end{equation}

\subsection{Canonical form} 
First we will identify the eigenvectors and the associated  eigenvectors 
of $A$ {\it formally} relying on (\re{unsXrw3}). Afterwords, we will 
prove the identifications rigorously. 
 
The set $\cW:=\{X\in\cV: \Lam X\in D(H)\}$ is dense 
in~$\cV$ under our conditions (1.3) and (1.4). 
Let us apply the both sides of identities~(\ref{unsXrw3}) 
to an arbitrary $X\in \mathcal W$. Using (\ref{unsF}), we formally obtain 
\beqn\label{unsX3} 
X&=&~~\int_{\mathbb R} dF(\omega)X=X_\varepsilon+X_0+X_a,
\nonumber\\
\nonumber\\
 AX&=&\!\!\!\!\!-i\int_{\mathbb R} \omega dF(\omega)X=AX_\varepsilon+AX_0+AX_a, 
\eeqn 
where 
\beqn 
X_\varepsilon &:=& \int_{|\omega|\ge \varepsilon} \Big[\Lambda_+^{-1}+ 
\frac{iP}{\omega}\Big] dE(\omega)\Lambda X
\nonumber\\
\nonumber\\
AX_\varepsilon &=& -i\int_{|\omega|\ge \varepsilon} \Big[\omega\Lambda_+^{-1}+ i P \Big] 
dE(\omega)\Lambda X\label{X}
\\
\nonumber\\
X_0 &:=& \Big[\Pi_\mathcal K- i P \int_{|\omega|\ge \varepsilon}\frac {dE(\omega)} {\omega}  \Lambda \Big]X,~~~~~~~~ AX_0=0,\label{X0} 
\\
\nonumber\\
X_a &:=& \Lambda_+^{-1}\Pi_0\Lambda X,  ~~~~~~~~~~~~~~~~~~~~~~~~~~~~~~~~~~~~~~
AX_a=P\Pi_0\Lambda X.~~~~~~~~~~~~~\label{Xa} 
\eeqn 
Here, (\ref{X}) means the expansion over the eigenvectors with eigenvalues 
$-i\omega$, while  (\ref{X0}), with the zero eigenvalue. 
Formula \eqref{Xa} means that 
$X_a$ is the associated eigenvector to the eigenvector 
$P \Pi_0\Lambda  X$, which corresponds to the zero eigenvalue. 
We justify the formal calculations  \eqref{unsX3}--\eqref{Xa} in the following lemma. 
 
\bl\label{lXXX} 
Formulas \eqref{unsX3}--\eqref{Xa} hold for $X\in\mathcal W$. 
\el 
\Pr 
i) Formulas \eqref{unsX3} and (\ref{X0}) are obvious. 
\medskip\\ 
ii) The last formula of (\re{Xa}) 
follows from the fact that $AX_a=J\Lam\Pi_0\Lam X=\Pi_\cK J\Lam\Pi_0\Lam X$, 
since  $\Lam J\Lam\Pi_0\Lam X=-iH\Pi_0\Lam X=0$ by definition of $\Pi_0$. 
\medskip\\ 
iii) 
Finally, let us prove (\re{X}). 
The representation (\re{unsXr}) implies that 
$\dot X(\cdot)\in C(\R,\cV)$ for $X(0)\in \cW$ since 
\be\la{dA} 
\dot X(t)=-i\Lam_+^{-1}e^{-iHt}H\Lam X(0)  + 
P \Pi_0\Lam X(0)+ 
P e^{-iH_Rt} \Pi_R \Lam X(0) 
\ee 
by Hille--Yosida's theorem \ci[Theorem 13.35 (c)]{Rudin}. 
On the other hand, according to (\re{ini}), 
\be\la{Xtd} 
\fr{X(t+\De t)-X(t)}{\De t}=A\fr{\ds\int_t^{t+\De t} X(s)ds}{\De t}~. 
\ee 
Here the left hand side converges to $\dot X(t)$ in $\cV$ 
as $\De t\to 0$,  since $\dot X(\cdot)\in C(\R,\cV)$, 
and the quotient on the right converges to 
$X(t)$ in $\cV$ by (\re{YD}). Hence, making $\De t\to 0$ we obtain 
\be\la{Xtd2} 
\dot X(t)=A X(t), 
\ee 
since 
the operator $A=JB$ is closed in $\cX$. 
Setting $t=0$ in (\re{dA}) and (\re{Xtd2}), and writing $X$ instead of $X(0)$, 
  we obtain 
\beqn\la{dA2} 
AX\!\!\!\!\!\!&=&\!\!\!\!\!-i\Lam_+^{-1}H\Lam X+P\Pi_0 \Lam X  + 
P \Pi_R\Lam X 
\nonumber\\ 
\nonumber\\ 
\!\!\!\!\!\!&=&\!\!\!\!\! 
-i\int_{|\om|\ge \ve} 
\Big[\om\Lam_+^{-1}+ 
i 
P \Big] 
dE(\om)\Lam X+ 
P \Pi_0 \Lam X,~~X\in\cW. 
\eeqn 
On the other hand, $AX=AX_\ve+AX_0+AX_a=AX_\ve+P \Pi_0 \Lam X$ by (\re{X0}) 
and (\re{Xa}). 
Hence,  (\re{dA2}) implies (\re{X}).\bo 
 
\bc\label{cJ2} 
The nontrivial Jordan blocks occur only for $\lambda=0$; they 
are of size $2\times 2$ {\rm (}in accordance with 
{\rm \cite[Proposition 5.1]{L1981}}{\rm )}, 
and their number is $\dimension \Ker H_R=\dimension [\Ker H\cap \mathcal R]$
where $H_R:=H|_{R\cap D(H)}$. 
This number is finite by~{\rm (1.4)}. 
\ec 
Further, we set $\Pi_\mathcal R:=1-\Pi_\mathcal K$ and introduce the `Green operator' 
\begin{equation}\label{G} 
G:=\Lambda_+^{-1}\Pi_\mathcal R+iP H_R^{-1}\Pi_R. 
\end{equation} 
It is continuous from $\mathcal X$ to $\mathcal V$ by our conditions 
(1.3) and {\rm (2.12)} 
according to Lemma \ref{lP}. Therefore, formulas  (\ref{X}) 
can be rewritten as 
\begin{equation}\label{XGZ} 
X_\varepsilon= 
G 
\int_{|\omega|\ge \varepsilon} 
dE(\omega)\Lambda X,\qquad 
AX_\varepsilon 
= 
-iG\int_{|\omega|\ge \varepsilon} 
\omega 
dE(\omega)\Lambda X,\qquad X\in\mathcal W, 
\end{equation} 
since the both integrals 
converge in $\mathcal X$ and 
belong to $R\subset\mathcal R$. 
 
\bc\label{ca} 
Let $h_k\in\mathcal X$ be an eigenfunction of $H$ corresponding to an 
eigenvalue $\omega_k\ne 0$. Then 
\begin{equation}\label{aGh} 
a_k:=Gh_k\in\mathcal V 
\end{equation} 
is the eigenfunction of $A$  corresponding to the 
eigenvalue $-i\omega_k$. 
 
\ec

\setcounter{equation}{0} 
\section{Application to eigenfunction expansion}

We are going to apply our results to justify the eigenfunction 
expansion (\ref{uns30i}) in the context of the system considered in \cite{KK11}. 
We have used this expansion for the calculation of the 
Fermi Golden Rule \cite[(5.14)]{KK11}. 
 
\subsection{Linearization at the kink} 
In \cite{KK11,KKm11} we studied the 
 1D relativistic Ginzburg--\allowbreak Landau equation 
\begin{equation}\label{GL} 
\ddot\psi(x,t)=\frac{d^2}{dx^2}\psi(x,t)+F(\psi(x,t)),\qquad x\in{\mathbb R} 
\end{equation} 
for real solutions $\psi(x,t)$. 
Here, $F(\psi)=-U'(\psi)$, where 
$U(\psi)$ is similar to the 
Ginzburg--\allowbreak Landau potential $U_{GL}(\psi) = (\psi^2-1)^2/4$, 
which corresponds to the cubic equation with $F(\psi) = \psi- \psi^3 $. 
Namely, $U(\psi)$ is a real smooth even function satisfying the following 
conditions: 
\be\label{U1} 
\left.
\ba{rcl}
U(\psi) &>&0,~~~~~~ \psi \ne \pm a\\ 
\\
U(\psi) &=& \frac{m^2}2(\psi\mp a)^2 + O(|\psi\mp a|^{14}),~~ x\to \pm a. 
\ea\right|
\ee 
where $a,m>0$. The main goal of \cite{KK11,KKm11} was to prove the asymptotic stability 
of solitons (kinks) $\psi(x,t)=s_v(x-vt)$ that 
move with constant velocity $|v|< 1$, and 
\begin{equation}\label{kink} 
s_v(x)\to\pm a,\qquad x\to\pm\infty. 
\end{equation} 
Substituting  $\psi(x,t)=s_v(x-vt)$ into (\ref{GL}), we obtain the corresponding stationary equation 
\begin{equation}\label{GLs} 
v^2 s_v''(x)=s_v''(x)+F(s_v(x)),\qquad x\in{\mathbb R}. 
\end{equation} 
The linearization 
of (\ref{GL}) 
at the kink $s_v(x-vt)$ in the moving frame 
reads as \eqref{NHl} with $X=(\psi,\dot\psi)\in L^2({\mathbb R})\otimes{\mathbb C}^2$ 
(for the corresponding complexification) 
and with the generator 
\cite[(4.6)]{KKm11} 
\begin{equation}\label{AA1} 
A_v=\left( 
\begin{array}{cc} 
v\frac d{dx}   & 1 \\ 
\frac{d^2}{dx^2}-m^2-V_v(x) & v\frac d{dx} 
\end{array} 
\right). 
\end{equation} 
Here, the potential 
\begin{equation}\label{Vv} 
V_v(x)=-F'(s_v(x))-m^2\in C^\infty({\mathbb R}). 
\end{equation} 
The kink $s_v(x)$ is an odd monotone function 
in a suitable coordinate $x$, 
while $F'(\psi)=-U''(\psi)$ is an even function of~$\psi$. 
Hence, the potential $V_v(x)$ is an even function of~$x$. 
Moreover, 
\begin{equation}\label{Vvm} 
|V_v(x)|\le Ce^{-\kappa|x|},\qquad x\in{\mathbb R}, 
\end{equation} 
where $\kappa>0$. 
The generator (\ref{AA1}) 
has the form $A_v=JB_v$ with 
\begin{equation}\label{BJ} 
B_v=\left( 
\begin{array}{cc} 
S_v & -v\frac d{dx} \\ 
  v\frac d{dx}   &   1 
\end{array} 
\right), 
\qquad
J:=\left( 
\begin{array}{cc} 
0 & 1 \\ 
-1 & 0 
\end{array} 
\right), 
\end{equation} 
where $S_v:=-\frac{d^2}{dx^2} +m^2+V_v(x)$.
Obviously, $JB_v\ne B_vJ$. 
Differentiating (\ref{GLs}), we obtain 
\begin{equation}\label{grst} 
\Bigl[S_v+v^2\frac{d^2}{dx^2}\Bigr]s_v'(x)=0. 
\end{equation}


 
\subsection{Spectral conditions} 
 
Conditions (1.4), {\rm (2.12)} hold for operators (\ref{BJ}) 
on $\mathcal X:=L^2({\mathbb R})\otimes{\mathbb C}^2$ 
by Lemma \ref{lcond23}. 
Condition (1.3)  for all $|v|<1$ 
follows   from Lemmas A.1 and A.2 
of \cite{KK13}. 
Here, we check (1.3) 
in the case $v=0$ for the completeness of the exposition. 
We will write $A$,  $B$ and $S$, respectively, instead of   $A_0$, $B_0$ and~$S_0$: 
\begin{equation}\label{AB0} 
A=\left( 
\begin{array}{cc} 
0   & 1 \\ 
-S & 0 
\end{array} 
\right), \qquad 
B=\left( 
\begin{array}{cc} 
S & 0 \\ 
0 & 1 
\end{array} 
\right),
\end{equation} 
where $ S:=-\frac{d^2}{dx^2}+m^2+ V_0(x)$.
The operators $B$ and $S$ are essentially selfadjoint 
in $L^2({\mathbb R})\otimes{\mathbb C}^2$ and 
$L^2({\mathbb R})$, respectively, by (\ref{Vvm}) and 
Theorems X.7 and X.8 of~\cite{RS2}. We will consider the closures of  $B$ and~$S$, 
which are both selfadjoint. In this case, 
\begin{equation}\label{LH0} 
\Lambda:=B^{1/2}=\left( 
\begin{array}{cc} 
\sqrt{S} & 0 \\ 
  0   &   1 
\end{array} 
\right),\qquad H:=\Lambda iJ\Lambda= 
i\left( 
\begin{array}{cc} 
0 & \sqrt{S} \\ 
 -\sqrt{S}  & 0 
\end{array} 
\right)=iJ\sqrt{S}. 
\end{equation} 
Hence, the operator $H$ is also selfadjoint on the domain $D(\sqrt{S})\oplus D(\sqrt{S})$. 
Thus, condition {\rm (2.21)} holds in our case.

 
\bl\label{lSB} 
Condition {\rm (1.3)} 
 holds 
for the operator $B$ on $\mathcal X=L^2({\mathbb R})\otimes{\mathbb C}^2$. 
\el 
\Pr 
Equation (\ref{grst}) with $v=0$ means that $\lambda=0\in \sigma_{pp}(S)$. 
Moreover, $\lambda=0$ is the minimal eigenvalue 
of $S$, since the corresponding eigenfunction $s_0'(x)$ 
does not vanish \cite[(1.9)]{KKm11}. Hence, 
\begin{equation}\label{sB} 
\sigma(S)\subset [0,\infty),\qquad \Ker S=(s_0'(x)). 
\end{equation} 
Further, 
the continuous spectrum of $S$ lies in $[m^2,\infty)$, and hence 
(\ref{sB}) implies 
\begin{equation}\label{spB0} 
\sigma(S)= 
 \{\lambda_0,\dotsc,\lambda_N\}\cup [m^2,\infty), 
\end{equation} 
where $0=\lambda_0<\dotsc<\lambda_N<m^2$. Finally, $\sigma(B)=\sigma(S)\cup \{1\}$, by~(1.3). 
\bo 
\medskip 
 
We will assume below the following spectral condition 
(imposed in \cite{KK11})) 
at the 
 edge point of the continuous spectrum of $S$: 
\begin{equation}\label{SC11} 
\mbox{\it The point $m^2$ is neither an eigenvalue 
nor a~resonance of $S$.} 
\end{equation} 
This  condition provides a regularity of the eigenvalue expansion 
(\re{uns30i}) at the edge points  $\pm m$ of the continuous spectrum.

\setcounter{equation}{0} 
\section{Orthogonal eigenfunction expansion} 
We are going to apply Proposition \ref{pX} 
to the case of operators (\ref{AB0}). First, 
(\ref{spB0}) implies that 
\begin{equation}\label{spB} 
\sigma(H) 
=(-\infty,-m]\cup 
 \{\omega_{-N},\dotsc,\omega_{-1}, \omega_0,\omega_1,\dotsc,\omega_N\} 
\cup [m,\infty), 
\end{equation} 
where $\omega_{\pm k}^2=\lambda_k$ for $k=0,\dotsc,N$. 
We denote by $\sigma_c=(-\infty, -m]\cup[m,\infty)$ the 
continuous spectrum of $H$, and 
\begin{equation}\label{PP} 
\Psi_0=\Lambda_+^{-1}\Pi_0\Lambda X(0),\qquad \Phi_0=P\Pi_0\Lambda X(0). 
\end{equation} 
Then $A\Psi_0=0$ and $A\Psi_0=\Phi_0$ 
by \eqref{Xa}, and 
hence, 
$t\Phi_0+\Psi_0$ is the solution to \eqref{NHl}. 
Now formula (\ref{unsXrw}) 
 can be rewritten as 
\beqn\label{uns302} 
X(t)=e^{At}X(0)&=&t\Phi_0+\Psi_0+ 
\sum_{-N}^N e^{-i\omega_k t} C_k a_k 
\nonumber\\
\nonumber\\
&&+ 
\int_{\sigma_c} \Bigl [\Lambda_\mathcal R^{-1}
+ 
\frac{iP}\omega\Bigr] e^{-i\omega t} d E(\omega)\Lambda X(0). 
\eeqn 
Here, $a_0\in\mathcal K$ and 
\begin{equation}\label{ak} 
a_k= 
\Bigl [\Lambda_\mathcal R^{-1}+\frac{iP}{\omega_k}\Bigr] h_k=Gh_k 
\in\mathcal X,\quad k\ne 0, 
\end{equation} 
where 
$h_k\in\mathcal R$ 
are the eigenfunctions 
of $H$ corresponding to the eigenvalues 
$\omega_k\ne 0$. By Corollary \ref{ca}, 
$a_k$ are the  eigenfunctions of $A$ 
corresponding to the eigenvalues 
$-i\omega_k$. 
\medskip 
 
Let us denote by $X^c(t)$ the integral 
in (\ref{uns302}): 
\begin{equation}\label{eife0} 
X^c(t) 
= 
\int_{\sigma_c}\Bigl [\Lambda_\mathcal R^{-1}+ 
\frac{iP}\omega\Bigr ] e^{-i\omega t} d E(\omega)\Lambda X(0). 
\end{equation} 
To prove (\ref{uns30i}), it remains to justify the eigenfunction expansion 
\begin{equation}\label{eife} 
X^c(t) 
= 
\int_{\sigma_c} e^{-i\omega t} C(\omega)\5 a_\omega\5 \5 d\omega, 
\end{equation} 
where $a_\omega$ are the generalized eigenfunctions of $A$ 
corresponding to the eigenvalues $-i\omega$,
and the meaning of the 
convergence of the integral will be specified later. 
Then 
(\ref{uns30i}) will follow from (\ref{uns302}). 
\medskip

By  (\ref{uns302}),  the function $X^c(t)$ 
is the solution to \eqref{NHl}, and hence 
\begin{equation}\label{Zc} 
Z^c(t):=\Lambda X^c(t)= 
\int_{\sigma_c} e^{-i\omega t} 
dE(\omega)\Lambda X(0) 
\end{equation} 
 is the solution 
to \eqref{CPF3}. 
We will deduce (\ref{eife}) from 
 the corresponding representation 
\begin{equation}\label{eife2} 
Z^c(t) 
= 
\int_{\sigma_c} e^{-i\omega t} C(\omega)\5 h_\omega\5 \5 d\omega, 
\end{equation} 
where $h_\omega$ are  generalized eigenfunctions of~$H$ corresponding to 
the eigenvalues $\omega$ normalized by
\begin{equation}\label{norm2} 
\langle h_{\omega}, h_{\omega'}\rangle=2\pi\5\delta(\omega-\omega'), 
\qquad \omega,\omega' \in\sigma_c. 
\end{equation} 
 The normalization  means by definition, that
\beqn\label{nor22} 
\left.
\ba{rcl}
\langle Z_1,Z_2 \rangle&=&2\pi 
\ds\int_{m\le  |\omega|\le M} C_1(\omega)\overline{C_2(\omega)}d\omega
\\
\\
\mbox{\rm for}~~~Z_k&=&\ds
\int_{m\le  |\omega|\le M}C_k(\omega) h_\omega d\omega\in\mathcal X
\ea\right|  
\eeqn 
For $\rho\in{\mathbb R}$ we denote
by  $L^2_{\rho}=L^2_{\rho}({\mathbb R})$  the 
weighted Hilbert space with the norm 
\begin{equation}\label{whs} 
\Vert\psi\Vert^2_{L^2_{\rho}}:=\int \langle x \rangle^{2\rho}|\psi(x)|^2dx,\qquad 
\langle x \rangle:=(1+x^2)^{1/2~}. 
\end{equation}

\bt\label{lor} Let condition \eqref{SC11} hold and $s>1$. Then, for 
$\omega\in \sigma_c$, 
there exists $h_\omega\in L^2_{-s}\otimes {\mathbb C}^2$  such that{\rm :} 
\medskip\\ 
i) 
$h_\omega$ is a continuous function of $\omega\in\sigma_c$ with values in 
 $L^2_{-s}\otimes {\mathbb C}^2$. 
\medskip\\ 
ii) The normalization (\re{norm2}) holds.
\medskip\\
iii) $h_\omega$ are the generalized eigenfunctions of $H$, i.e., 
\begin{equation}\label{eigf2} 
H Z=\int_{\sigma_c}  \omega\5 
C(\omega)\5 h_\omega\5 \5 d\omega \qquad \text{if} \quad Z= \int_{\sigma_c} 
C(\omega)\5 h_\omega\5 \5 d\omega 
\in D(H). 
\end{equation} 
iv) 
The eigenfunction expansion \eqref{eife2} 
holds in the following sense{\rm :} 
\begin{equation}\label{icon} 
\Bigl\Vert Z^c(t)-\int_{m\le  |\omega|\le M} e^{-i\omega t} C(\omega)\5 h_\omega 
\5 \5 d\omega\Bigr\Vert_{L^2\otimes{\mathbb C}^2}\to 0, \qquad M\to\infty. 
\end{equation} 
\et 

\Pr 
i)
We construct the generalized eigenfunctions and 
the eigenfunction expansion \eqref{eife2} 
by solving equation (\re{CPF3})
for $Z^c(t)=(Z_1^c(t),Z_2^c(t))$.
By (\ref{LH0}), the
equation  is equivalent 
to 
the system 
\begin{equation}\label{Z2} 
\dot Z_1^c(t)=\sqrt{S}Z_2^c(t), \qquad 
\dot Z_2^c(t)=-\sqrt{S}Z_1^c(t). 
\end{equation} 
Eliminating $Z_2^c(t)$, we obtain 
\begin{equation}\label{Z3} 
\ddot Z_1^c(t)= -SZ_1^c(t). 
\end{equation} 
Further
we apply Theorem XI.41 of \cite{RS3} and
the arguments of \cite[pp 114-115]{RS3}. Namely,
the rapid decay (\re{Vvm}) and our spectral condition (\re{SC11}) imply
the following Limiting Absorption Principle (LAP) \cite{A,KKW12,RS3}:
\begin{equation}\la{LAP}
R(\lambda\pm i\varepsilon)\to R_\pm(\lambda),\qquad \varepsilon\to+0,
~~~~\lambda\in [m^2,\infty),
\end{equation}
where
$R(z):=(S-z)^{-1}$ and
the convergence holds
in the strong topology
of the space of continuous operators
$L^2_{s}\to L^2_{-s}$ with $s>1$.
Moreover, the traces of the resolvent $R_\pm(\lambda)$ are continuous functions of
$\lambda\ge m^2$
with values in $L(L^2_{s}, L^2_{-s})$.
The continuity at $\lambda>0$ has been 
established by Agmon,  see \cite{A,KKW12}.
The continuity at $\lambda=0$ under condition (\re{SC11}) is proved in 
\ci[formulas (3.12)]{Kop2011}.
The LAP serves as the basis for
the eigenfunction expansion 
\beqn\la{eife3}
Z_1^c(t)&=&\int_{\si_c}
d \cE(\omega^2)[Z_1^c(0)\cos \omega t+Z_2^c(0)\sin \omega t] 
\nonumber\\
\nonumber\\
&=&
\int_{\si_c} e^{-i\omega t} C(\omega)\5 e_\omega\5 \5 d\omega,
\eeqn
where  $d\cE(\lambda)$ is the spectral resolution of $S$,
while
$e_\omega\in L^2_{-s}$
are generalized eigenfunctions
of $S$ corresponding to the eigenvalues $\omega^2\ge m^2$.
Here the first identity follows by Spectral Theorem, while the second
follows 
by  Theorem XI.41 (e) of \cite{RS3}.
The eigenfunctions are defined by formulas
of \cite[pp 114-115]{RS3}:
\begin{equation}\la{eigf}
e_\omega=W^*(\omega)f_\omega,~
~~ W(\omega):=[1+VR_0(\omega^2+i0)]^{-1},
~ \omega\in \sigma_c.
\end{equation}
where 
$f_\omega(x):=\sin|\omega| x$ and
$R_0(\lambda):=(-\De+m^2-\lambda)^{-1}$.

The operator $W(\omega)$
 is a~continuous function of $\omega\in\sigma_c$
with values in
 $L(L^2_{s},L^2_{s})$ by the formula 
\be\la{LL2}
[1+VR_0(\lambda)]^{-1}=1-VR(\lambda)
\ee
and the decay (\re{Vvm}).
Respectively,
the adjoint operator $W^*(\omega)$
is a continuous function of $\omega\in \sigma_c$
with values in
 $L(L^2_{-s},L^2_{-s})$.
As the result, $e_\omega$ is a
continuous function of $\omega\in\sigma_c$
with values in
$L^2_{-s}$. The normalization of $e_\omega$ coincides with the same of
the 'free' generalized eigenfunctions
$f_\omega$:
\begin{equation}\la{norm}
\langle e_{\omega}~, e_{\omega'}\rangle=\pi\5\delta(|\omega|-|\omega'|)~,
\quad \omega,\omega' \in\sigma_c,
\end{equation}
which follows from the last formula on page 115 of \cite{RS3}.
Finally,
Theorem XI.41 (e) of
\cite{RS3} implies that
the last integral (\re{eife3}) converges in
$L^2=L^2(\R)$:
\begin{equation}\la{icon0}
\Vert Z_1^c(t)-\int_{m\le |\omega|\le M} e^{-i\omega t} C(\omega)
\5 e_\omega\5 
\5 d\omega\Vert_
{L^2}
\to 0~,\qquad M\to\infty.
\end{equation}

Now (\re{eife2}) for $Z_1^c(t)$ follows from (\re{eife3}).
For $Z_2^c(t)$ we  use the first
equation of (\re{Z2}),  which implies
\beqn\la{eife4}
Z_2^c(t)
=
-i
\int_{\sigma_c} \sgn\omega \5\5\5e^{-i\omega t} C(\omega)\5 e_\omega\5 \5 d\omega.
\eeqn
Combining (\re{eife3}) and (\re{eife4}), we obtain (\re{eife2})
with
\begin{equation}\la{Zk}
h_\omega:=\left(\ba{c} 1\\-i ~
\sgn\omega\ea\right)e_\omega,
\end{equation}
which is the continuous function of $\om\in\si_c$ with values in 
 $L^2_{-s}\otimes {\mathbb C}^2$.
\medskip\\
ii) Normalization (\re{norm2}) follows from (\re{norm}).
\medskip\\
iii) 
 $Z^c(t)\in D(H)$ means that $Z^c_{1,2}(t)\in D(\sqrt{S})$. 
Furthermore,
\begin{equation}\la{eigf3}
H Z^c(t)=i\sqrt{S}
\left(\ba{c}
 Z^c_2(t)\\ -Z^c_1(t)
\ea
\right).
\end{equation}
Now (\re{eigf2}) follows from the expansions (\re{eife3}) and (\re{eife4})
for $Z^c_{1,2}(t)$ by
\cite[Theorem XI.41 (c)]{RS3}, since $e_\omega$ are the generalized 
eigenfunctions of $S$
with the eigenvalues $\omega^2$, and {\it formally},
\begin{equation}\la{eigf4}
i\sqrt{S}
\left(\ba{c}
 -i~\sgn\omega\\ -1
\ea
\right)e_\omega=
\left(\ba{c}
\sgn\omega\\ -i
\ea
\right)|\omega|e_\omega=\omega h_\omega~.
\end{equation}
iv) (\re{icon}) follows from (\re{icon0}) and similar convergence for
$Z_2^c$.\bo

\setcounter{equation}{0} 
\section{Nonorthogonal eigenfunction expansion} 
Let us denote by  $Z^c_M(t,x)$ the integral in (\ref{icon}). 
This integral is defined for almost all~$x$; i.e., 
\begin{equation}\label{icon4} 
 Z^c_M(t,x):= 
\int_{m\le  |\omega|\le M} e^{-i\omega t} C(\omega)\5 h_\omega(x) 
\5 \5 d\omega, \qquad \text{a.a.}\quad x\in{\mathbb R}. 
\end{equation} 
To justify (\ref{eife}) we should adjust the meaning of this integral 
relying on 
 the following lemma, which is proved in~\cite{KK13}. 
 
\bl [{\rm \cite[Lemma 5.1]{KK13}}] \label{lsi} 
Let condition {\rm \eqref{SC11}} hold. 
Then 
\medskip\\ 
i) 
The integral {\rm (\ref{icon4})} 
converges absolutely in $L^2_{-s}\otimes{\mathbb C}^2$ for every $s >1$: 
\begin{equation}\label{whs2} 
\int_{m\le  |\omega|\le M} \Vert C(\omega)\5 h_\omega\Vert_{L^2_{-s}\otimes {\mathbb C}^2}\5 d\omega 
<\infty, \qquad M>m. 
\end{equation} 
ii) The integral   of these 
$L^2_{-s}\otimes{\mathbb C}^2$-valued functions over $m\le  |\omega|\le M$ 
coincides with 
{\rm (\ref{icon4})} almost everywhere. 
\el 
 
Further, we express $X^c(t)$ in terms of $Z^c(t)$ and the Green operator (\ref{G}), 
and prove the appropriate continuity of $G$, which allows us to deduce 
 (\ref{eife}) from  (\ref{icon}).

\subsection{Reconstruction via the Green operator} 
Similarly to (\ref{XGZ}), we use  (\ref{Zc}) to rewrite integral (\ref{eife0}) as 
\begin{equation}\label{X2} 
X^c(t)= 
G 
\int_{\sigma_c} 
e^{-i \omega t} 
dE(\omega)\Lambda X(0)=G Z^c(t), 
\end{equation} 
taking into account that $Z^c(t)\in R\subset\mathcal R$ 
and that the Green operator 
$G:\mathcal X\to\mathcal V$ is continuous. 
Now  \eqref{eife2} implies that 
\begin{equation}\label{X3} 
X^c(t)= 
G 
\int_{\sigma_c} 
e^{-i \omega t} 
C(\omega) h_\omega d\omega. 
\end{equation} 
Similarly to (\ref{XGZ}), 
\beqn\label{X4} 
AX^c(t)&\!\!\!\!\!\!=&\!\!\!\!\!\!-i 
G 
\int_{\sigma_c} 
e^{-i \omega t} 
\omega\5 dE(\omega)\Lambda X(0) 
\nonumber\\ 
\nonumber\\ 
&\!\!\!\!\!\!=&\!\!\!\!\!\!-i 
G 
H\int_{\sigma_c} 
e^{-i \omega t} 
 dE(\omega)\Lambda X(0)=-i 
G 
H Z^c(t),~~  X(0)\in\mathcal W. ~~  ~~ ~~~~ ~~
\eeqn 
Therefore,  (\ref{eigf2}) 
gives 
\begin{equation}\label{X5} 
AX^c(t)=-i 
G 
\int_{\sigma_c} 
e^{-i \omega t} 
\omega\5 C(\omega) h_\omega d\omega. 
\end{equation} 
\subsection{Continuity of the Green operator} 
Now we are going to establish the continuity of the Green operator 
$G$ in the weighted norms (\ref{whs}). 
We will simplify the form of  $G$ in the concrete case (\ref{LH0}) by 
proving that 
\begin{equation}\label{PHR} 
P H_R^{-1}=\Pi_{\mathcal K}J\Lambda H_R^{-1}=0. 
\end{equation} 
First we note that 
$$ 
\Ker H=\Ker S\oplus \Ker S,\quad R=\Ran H=\Ran S\oplus\Ran S 
$$ 
by (\ref{LH0}). Further, we set 
$$ 
S_+:=S|_{\Ran S}:\Ran S\cap D(S)\to \Ran S, 
$$ 
and let 
$P_0$ denote the orthogonal projection of $L^2({\mathbb R})$ onto 
$\Ker S$. Then $P_+:=1-P_0$ is the orthogonal projection of $L^2({\mathbb R})$ onto 
$\Ran S$, and now (\ref{LH0}) implies 
$$
\ba{ll}
\Pi_{\mathcal K}=\left(\begin{array}{cc} 
P_0 & 0 \\ 
  0   &  0 
\end{array} 
\right),&
\Pi_{\mathcal R}=\left( 
\begin{array}{cc} 
P_+ & 0  \\ 
 0 &  1 
\end{array} 
\right),\quad 
\\
\\
\Pi_0=\left( 
\begin{array}{cc} 
P_0 & 0  \\ 
0&  P_0 
\end{array} 
\right),&
\Pi_R=\left( 
\begin{array}{cc} 
P_+ & 0  \\ 
0&  P_+ 
\end{array} 
\right). 
\ea
$$
Hence, finally, 
\begin{equation}\label{HRi} 
H_R^{-1}=i\left( 
\begin{array}{cc} 
0 & S_+^{-1/2}  \\ 
 - S_+^{-1/2}&  0 
\end{array} 
\right) 
\end{equation} 
by (\ref{LH0}), and therefore, 
\beqn\label{JLH} 
J\Lambda H_R^{-1} &=& 
i\left(\begin{array}{cc} 
0 & 1 \\ 
-1   &  0 
\end{array} 
\right) 
\left(\begin{array}{cc} 
S_+^{1/2} & 0 \\ 
  0   &   1 
\end{array} 
\right) 
\left( 
\begin{array}{cc} 
0 & S_+^{-1/2}  \\ 
 - S_+^{-1/2}&  0 
\end{array} 
\right) 
\nonumber\\
\nonumber\\
&=& 
-i\left(\begin{array}{cc} 
S_+^{-1/2} & 0 \\ 
  0   &   P_+ 
\end{array} 
\right). 
\eeqn 
Applying $\Pi_\cK$, we get  (\ref{PHR}). 
\medskip

Now definition (\ref{G})  implies that 
\begin{equation}\label{G1} 
G=\Lambda^{-1}_+\Pi_\mathcal R. 
\end{equation} 
The following lemma is a generalization 
of \cite[Lemma 5.2]{KK13}. 
 
\bl\label{lambda} 
The operator $G:L^2_{\rho}\otimes{\mathbb C}^2\to L^2_{\rho}\otimes{\mathbb C}^2$ 
is continuous for every $\rho\in{\mathbb R}$. 
 
\el 
\Pr 
Using the first formula of (\ref{LH0}) and the  formula for $\Pi_\cR$, we get 
\begin{equation}\label{pi0} 
\Lambda^{-1}_\mathcal R\Pi_\mathcal R= 
\left( 
\begin{array}{cc} 
S_+^{-1/2}P_+ & 0 \\ 
  0   &   1 
\end{array} 
\right) 
=\left( 
\begin{array}{cc} 
QP_+ & 0 \\ 
  0   &   1 
\end{array} 
\right),  
\end{equation} 
where $Q:=(SP_++P_0)^{-1/2}$.
Hence, it suffices to prove the continuity of 
the operator $Q_+P_+$ in $L^2_{\rho}$, which means the 
continuity of operator 
\begin{equation}\label{QPDO} 
\langle x\rangle^{\rho} Q P_+ 
\langle x\rangle^{-\rho}:~~ L^2({\mathbb R})\to L^2({\mathbb R}). 
\end{equation} 
To prove this continuity, we note that 
$Q$ is a~PDO of the class $HG^{-1,-1}_1$, see Definition 25.2 in \cite{Shubin}. 
This fact follows from 
\cite[Theorem 29.1.9]{Hor4} 
and also by 
an extension of \cite[Theorem 11.2]{Shubin} to PDOs 
with nonempty continuous spectrum. It is important that 
operator $Q$ is a~PDO with the main symbol $\xi^2$, and 
$$ 
\xi^2\not\in (-\infty,0],\ \ \ \xi\ne 0;\qquad \sigma(S_+)\cap (-\infty,0]=\emptyset 
$$ 
by (\ref{sB}). 
Hence, conditions (10.1) and (10.2) 
of \cite{Shubin} hold. 
 
Now the continuity (\ref{QPDO}) 
follows by the Theorem of Composition of the PDO. 
\bo 
\medskip 
 
Lemma \ref{lambda} with $\rho=-s$ and Lemma \ref{lor} 
imply (cf. (\ref{ak})) that 
\begin{equation}\label{Xom} 
a_\omega:= 
G 
h_\omega\in L^2_{-s}\otimes{\mathbb C}^2,\qquad  s>1. 
\end{equation} 
Now we can prove the following lemma. 
 
\bl\label{leig} 
$a_\omega$ are the generalized eigenfunctions of $A$ corresponding to 
the eigenvalues $-i\omega$. 
\el 
\Pr 
Formulas (\ref{X3}) and (\ref{X5}) imply that 
\begin{equation}\label{X6} 
X^c(t)= 
\int_{\sigma_c} 
e^{-i \omega t} 
C(\omega) a_\omega d\omega, 
~~
AX^c(t)= 
\int_{\sigma_c} 
e^{-i \omega t} 
\omega\5 C(\omega) a_\omega d\omega, 
\end{equation} 
for $X^c(0)\in\mathcal W $
by definition (\ref{Xom}), 
Lemma \ref{lsi} and the last corollary with $\rho=-s<-1$. 
These identities mean that 
$a_\omega$ are the generalized eigenfunctions 
 in the sense of \cite[(80b)]{RS3}. 
\bo 
\medskip 
 
Finally, the main result of our paper is the following. 
 
\bt\label{tmain} 
Let condition {\rm \eqref{SC11}} hold, $X(0)\in \mathcal V$ and $s>1$. Then the eigenfunction 
expansion \eqref{eife} holds in the following sense {\rm (}cf.\ \eqref{icon}{\rm )}: 
\begin{equation}\label{icon2} 
\Bigl\Vert X^c(t) 
-\int_{m\le  |\omega|\le M} e^{-i\omega t} C(\omega)\5 a_\omega\5 \5 d\omega\Bigr\Vert_\mathcal V\to 0,\qquad M\to\infty, 
\end{equation} 
where the integral converges 
in  $L^2_{-s}\otimes{\mathbb C}^2$, and hence a.e. as in  {\rm (\ref{icon4})}.

\et 
\Pr Formulas (\ref{X2}) and (\ref{Xom}) imply that 
\beqn\label{icon22} 
&&X^c(t) 
-\int_{m\le  |\omega|\le M} e^{-i\omega t} C(\omega)\5 a_\omega\5 \5 d\omega 
\nonumber\\
\nonumber\\
&=&G\Big[ 
Z^c(t) 
-\int_{m\le  |\omega|\le M} e^{-i\omega t} C(\omega)\5 h_\omega\5 \5 d\omega 
\Big]. 
\eeqn 
Therefore, (\ref{icon2}) follows from (\ref{icon}), because the Green operator 
$G:\mathcal X\to\mathcal V$ is continuous. 
\bo

\setcounter{equation}{0} 
\section{Symplectic normalization} 
Now let us renormalize $h_\omega$ as follows: 
\begin{equation}\label{nor} 
\langle h_{\omega}, h_{\omega'}\rangle=|\omega|\5\delta(\omega-\omega'), \qquad \omega,\omega' \in \sigma_c. 
\end{equation} 
This means that 
\be\label{nor2} 
\left.
\ba{rcl}
\langle Z_1,Z_2 \rangle&=&
\ds\int_{m\le |\omega|\le M}|\omega|\5 C_1(\omega)\overline{C_2(\omega)}d\omega  \\
\\
\mbox{for}~~ Z_k&=&\ds\int_{m\le |\omega|\le M}C_k(\omega) h_\omega d\omega\in\mathcal X
\ea
\right|
\ee
similarly to (\re{nor22}).
We will express these formulas 
in terms of $X_k:=G Z_k\in\mathcal V$ 
and the eigenfunctions $a_\omega:=G h_\omega$. First, 
\begin{equation}\label{nor3} 
X_k=\int_{m\le |\omega|\le M}C_k(\omega) a_\omega d\omega 
\end{equation} 
by Lemma \ref{lambda}. 
Further, $Z_k\in R$, and so (\ref{eigf2}),  \eqref{nor2}  imply that 
\begin{equation}\label{nor4} 
\langle H_R^{-1}Z_1,Z_2 \rangle=\int_{m\le |\omega|\le M}\sgn\omega\5\5 C_1(\omega)\overline{C_2(\omega)}d\omega. 
\end{equation} 
On the other hand, this scalar product can be expressed in $X_k$. 
\bl 
Let $Z_1, Z_2$ be defined as in \eqref{nor2}. Then 
\begin{equation}\label{nor5} 
\langle H_R^{-1}Z_1,Z_2 \rangle=-i\langle X_1,J X_2 \rangle. 
\end{equation} 
\el 
\Pr 
First, $Z_1,Z_2\in R\subset \mathcal R$, and hence, 
$$ 
\Pi_{\mathcal R}Z_k=Z_k. 
$$ 
Now (\ref{G1}) implies (\ref{nor5}): 
$$ 
\ba{rcl}
\langle X_1,J X_2 \rangle&=& 
\langle GZ_1,JGZ_2\rangle= 
\langle \Lambda^{-1}_{\mathcal R}\Pi_{\mathcal R}Z_1,J \Lambda^{-1}_{\mathcal R}\Pi_{\mathcal R}Z_2\rangle 
\\
\\
&=&-\langle\Lambda^{-1}_{\mathcal R}\Pi_{\mathcal R} J\Lambda^{-1}_{\mathcal R}Z_1,Z_2\rangle 
=i\langle  H_R^{-1}Z_1,Z_2\rangle, 
\ea
$$ 
since 
$$ 
\ba{rcl}
\Lambda^{-1}_{\mathcal R}\Pi_{\mathcal R} J\Lambda^{-1}_{\mathcal R} 
&=&\left( 
\begin{array}{cc} 
S_+^{-1/2}P_+ & 0 \\ 
  0   &   1 
\end{array} 
\right) 
\left( 
\begin{array}{cc} 
0 & 1 \\ 
  -1   &   0 
\end{array} 
\right) 
\left( 
\begin{array}{cc} 
S_+^{-1/2} & 0 \\ 
  0   &   1 
\end{array} 
\right)
\\
\\
&=&\left( 
\begin{array}{cc} 
0 & S_+^{-1/2} \\ 
-S_+^{-1/2}     &   0 
\end{array} 
\right)=-i H_R^{-1} 
\ea
$$ 
by  the first formula of (\ref{LH0}) and by (\ref{pi0}) and~(\ref{HRi}). 
\bo 
\medskip 
 
Using this lemma and  (\ref{nor4}), we get
\begin{equation}\label{nor20} 
-i\langle X_1,JX_2 \rangle= 
\int_{m\le |\omega|\le M}\sgn\omega\5\5 C_1(\omega)\overline{C_2(\omega)}d\omega 
. 
\end{equation} 
By definition, (\ref{nor3}) and (\ref{nor20}) mean that 
\begin{equation}\label{nor21} 
\langle a_{\omega},~ Ja_{\omega'}\rangle=i\sgn\omega\5\5\5 \delta(\omega-\omega'), 
\qquad \omega,\omega'\in \sigma_c. 
\end{equation} 
Now expansion \eqref{eife}  coincides with \cite[(2.1.13)]{BS03}, thereby justifying our calculation 
of the Fermi Golden Rule 
for all solutions without the antisymmetry condition imposed in \cite{KK11}.

\appendix

\setcounter{equation}{0} 
\protect\renewcommand{\thesection}{\Alph{section}} 
\protect\renewcommand{\theequation}{\thesection.\arabic{equation}} 
\protect\renewcommand{\thesubsection}{\thesection.\arabic{subsection}} 
\protect\renewcommand{\thetheorem}{\Alph{section}.\arabic{theorem}}

\section{Examples} 
 
Let us show that 
conditions 
(1.4), {\rm (2.12)} and {\rm (3.1)} hold for elliptic  PDO 
\begin{equation}\label{Sm} 
P\psi(x)=\int e^{-ix\xi}P(x,\xi)\hat\psi(\xi)d\xi, 
\end{equation} 
which are the main objects of the theory. 
We will use the classes $\mathcal S^m$ of PDO similar to the ones 
introduced in \cite{Grushin1970}. 
\bd 
i) $P\in\mathcal S^m$ if, for any multiindices $\alpha,\beta$, 
\begin{equation}\label{Sm2} 
\sup_{x\in{\mathbb R}^n}|(1+|x|)^N\partial_\xi^\alpha \partial_x^\beta P(x,\xi)|\le C_{\alpha\beta N} 
(1+|\xi|)^{m-|\alpha|}, \qquad \xi\in{\mathbb R}^n 
\end{equation} 
with $N=0$ for $\beta=0$ and any $N>0$ for $\beta\ne 0$. 
\medskip\\ 
ii)  $P\in\mathcal S^m_0$ if  \eqref{Sm2} holds for any  multiindices $\alpha,\beta$ and all $N>0$. 
\medskip\\ 
iii) $P\in \mathcal S^m$ is elliptic of order $m$ if $P=P_m+R$, where 
 $P_m\in\mathcal S^m$ and 
\begin{equation}\label{Sm3} 
|P_m(x,\xi)|\ge C(1+|\xi|)^m, \qquad x,\xi\in{\mathbb R}^n, 
\end{equation} 
while $R\in\mathcal S^\mu_0$ with $\mu<m$. 
\ed 
 
Let  $\mathcal H^s=\mathcal H^s({\mathbb R}^n)$ denote the Sobolev spaces, and 
$\mathcal X=L^2({\mathbb R}^n)$. 
 Any operator 
$P\in \mathcal S^m$ is continuous $H^s\to H^{s-m}$ for 
 $s\in{\mathbb R}$, see Theorem 3.1 of 
\cite{Grushin1970}.

\bl\label{lcond23} 
Let 
$B\in \mathcal S^m$  be an elliptic  PDO of order $m$ which is symmetric 
on $C_0^\infty({\mathbb R}^n)$, and 
let $J\in \mathcal S^0$ be an elliptic  PDO of order $0$ which is antisymmetric 
on $C_0^\infty({\mathbb R}^n)$. 
Then 
\medskip\\ 
i) $B$ {\rm (}respectively, $J${\rm )} is selfadjoint {\rm (}respectively, skew selfadjoint{\rm )} operator 
with domain 
\begin{equation}\label{Sm4} 
D(B)=\mathcal H^m,\qquad D(J)=\mathcal X. 
\end{equation} 
ii) Condition {\rm (1.4)}   holds. 
\medskip\\ 
iii) Condition {\rm (2.12)} holds. 
\medskip\\ 
iv) Condition {\rm (3.1)} holds. 
\el 
\Pr i) 
The  Fredholm theory of elliptic PDO   on ${\mathbb R}^n$ 
\cite[Section 25.4]{Shubin} implies that $B\psi\in\mathcal X$ 
if and only if $\psi\in \mathcal H^m$, 
and the same is true for $B^*$. Hence, $D(B^*)=D(B)$, and therefore, $B^*=B$. 
Similarly, $J^*=-J$. 
\medskip\\ 
ii) 
The  Fredholm theory of elliptic PDOs  on ${\mathbb R}^n$ 
 implies that the space 
$\mathcal K:=\Ker B$ is finite dimensional and 
$\mathcal K\subset \mathcal H^s$ for any $s\in{\mathbb R}$. 
Hence, (1.4) holds. 
\medskip\\ 
iii) 
The operator $B_+:=B+\Pi_\mathcal K$ and its 
main symbol $B_+^m(x,\xi)$ satisfy 
$$ 
B_+^m(x,\xi)\not\in (-\infty,0], \quad  \xi\ne 0;\qquad \sigma(B_+)\cap (-\infty,0]=\emptyset 
$$ 
by (1.3). 
Therefore, conditions (10.1) and (10.2) 
of \cite{Shubin} hold for $B_+$, and 
hence, 
$\Lambda_+:=\sqrt{B_+}\ge 0$ 
is also an elliptic PDO of class $\mathcal S^{m/2}$. 
This follows similarly to Theorem 29.1.9 of 
\cite{Hor4} 
and also by 
an extension of Theorem 11.2 of \cite{Shubin} to PDO 
with nonempty continuous spectrum. 
Finally, $\Lambda_+=\Lambda+\Pi_\mathcal K$. Therefore, $\mathcal V=\mathcal H^{m/2}$, 
and hence {\rm (2.12)} holds, inasmuch as $J\in \mathcal S^{0}$. 
\medskip\\ 
iv) The operator $H$ is elliptic PDO of class $\cS^{m}$ 
by the theorem of composition. 
It is  obviously symmetric on $C_0^\infty(\R^n)$, 
and hence $H$~is selfadjoint on the domain $\mathcal H^m$ by 
the argument above. Thus {\rm (3.1)} is established. 
\bo

\br\label{rpdo} 
i) An example of elliptic operators $B\in \mathcal S^2$ and $J\in\mathcal S^0$ 
satisfying all conditions \eqref{Hs}, {\rm (1.3), (1.4)}, {\rm (2.12)}, 
{\rm (3.1)} 
is provided in Lemma~{\rm \ref{lSB}}. 
\medskip\\ 
ii) 
In the framework of Lemma~{\rm \ref{lcond23}} we should take $m\ge 0$ 
to keep condition~{\rm (1.3)}. 
\medskip\\ 
iii) The last condition of  \eqref{Hs} implies that 
the order of $J$ should be zero. 
 
\er 
 

\end{document}